\documentclass[11pt]{amsart}
\author{Christian J. Berghoff}
\title{Elliptic Gau{\ss} sums and Schoof's algorithm}
\address{Universit\"at Bonn, Mathematisches Institut, Endenicher Allee 60, 53115 Bonn, Germany}
\email{berghoff@math.uni-bonn.de}

\usepackage{lmodern}
\usepackage[T1]{fontenc}
\usepackage[utf8x]{inputenc}
\usepackage{a4wide}
\usepackage[english]{babel}
\usepackage{amsmath, amssymb, amsthm, amsfonts}
\usepackage{framed}
\usepackage{enumerate}

\usepackage{mathabx}
\usepackage[all]{xy}
\usepackage{textcomp}
\usepackage{stmaryrd}
\usepackage{nicefrac}
\usepackage[numbers]{natbib}
\usepackage[pdftex]{graphicx}
\usepackage{layout}
\usepackage{fancyhdr}
\usepackage{setspace}
\usepackage{color}
\usepackage{hyperref}

\usepackage{float}
\usepackage{pgfplots}
\pgfplotsset{compat=1.8}
\usepackage{subfigure}

\numberwithin{equation}{section}

\theoremstyle{plain}

\newtheorem{satz}{Theorem}[section]
\newtheorem{prop}[satz]{Proposition}
\newtheorem{lemma}[satz]{Lemma}

\theoremstyle{definition}

\newtheorem{defi}[satz]{Definition}

\theoremstyle{remark}
\newtheorem{bem}[satz]{Remark}

\DeclareMathOperator{\gal}{Gal}

\DeclareMathOperator{\deg2}{deg}
\DeclareMathOperator{\ord}{ord}

\DeclareMathOperator{\M}{\mathsf{{M}}}
\DeclareMathOperator{\Cpoly}{\mathsf{C}}

\newcommand{\OO}{\mathcal{O}}

\newcommand{\F}{\mathbb{F}}
\newcommand{\Ln}{\mathbb{L}}
\newcommand{\N}{\mathbb{N}}
\newcommand{\Z}{\mathbb{Z}}

\newcommand{\Cn}{\mathbf{C}}
\newcommand{\Bn}{\mathbf{B}}
\newcommand{\An}{\mathbf{A}}

\newcommand\numberthis{\addtocounter{equation}{1}\tag{\theequation}}
\newcommand{\mmid}{\mid \mid}

\begin{document}

\begin{abstract}
We present a new approach to handling the case of Atkin primes in Schoof's algorithm for counting points on elliptic curves over finite fields. Our approach is based on the theory of polynomially cyclic algebras developed in \cite{MiVu}, which we recall as far as necessary, and was elaborated in \cite{ma}. We then proceed to describe our method, which essentially relies on transferring costly computations in extensions of $\mathbb{F}_p$ to isomorphic ones endowed with a special structure allowing to reduce run-time. We analyse the new run-time and conclude this procedure yields some improvement as compared to the classical approaches.
\end{abstract}
\maketitle

\tableofcontents
\section{Introduction}
In this document we consider a new approach for a building block of Schoof's algorithm, which computes the number of rational points of an elliptic curve $E$ over a finite field $\F_p$. In order to do so, we first recall some details on so-called polynomially cyclic algebras that were defined in \cite{MiVu} before proceeding to present our new method in section \ref{sec:anwendung}. A large proportion of the new content presented in this section stems from the master's thesis \cite{ma} following first considerations to this end effected in \cite{MiVu}. Whereas former improvements in run-time for Schoof's algorithm essentially result from reducing the degree of the extensions of $\F_p$ in which computations are performed, our overall strategy is to transfer calculations to an isomorphic extension the structure of which allows to make the step dominating the run-time much more efficient.\\

Within this work we will only consider primes $p>3$ and thus assume that the curve in question is given in the Weierstra{\ss} form
\[ 
E: Y^2=X^3+aX+b=f(X),
 \]
where $a, b \in \F_p$. We will always identify $E$ with its set of points $E(\overline{\F_p})$. For the following well-known statements cf. \cite{Silverman, Washington}. We assume that the elliptic curve is neither singular nor supersingular. It is a standard fact that $E$ is an abelian group with respect to point addition. Its neutral element, the point at infinity, will be denoted $\OO$. For a prime $\ell \neq p$, the $\ell$-torsion subgroup $E[\ell]$ has the shape
\[ 
E[\ell]\cong \Z/\ell\Z \times \Z/\ell\Z.
 \]
Using the addition formulae for $k \in \N$ one can derive polynomials $G_k, H_k \in \F_p[X]$ such that
\begin{equation}\label{eq:multi_poly}
k(X, Y)=(G_k(X), Y H_k(X)).
 \end{equation}
In the endomorphism ring of $E$ the Frobenius homomorphism
\[ 
\phi_p: (X, Y) \mapsto (\varphi_p(X), \varphi_p(Y))=(X^p, Y^p)
 \]
satisfies the quadratic equation 
\begin{equation}\label{eq:char_gl}
0=\chi(\phi_p)=\phi_p^2-t\phi_p+p,
 \end{equation} 
where $|t|\leq 2\sqrt{p}$ by the Hasse bound. By restriction $\phi_p$ acts as a linear map on $E[\ell]$. The number of points on $E$ over $\F_p$ is given by $\#E(\F_p)=p+1-t$ and is thus immediate from the value of $t$.\\
The idea of Schoof's algorithm now consists in computing the value of $t \mod \ell$ for sufficiently many small primes $\ell$ by considering $\chi(\phi_p) \mod \ell$ and in afterwards combining the results by means of the Chinese Remainder Theorem. In the original version this requires computations in extensions of degree $O(\ell^2)$.\\
However, a lot of work has been put into elaborating improvements. Let $\Delta=t^2-4p$ denote the discriminant of equation \eqref{eq:char_gl}. Then we distinguish the following cases:
\begin{enumerate}
\item If $\left( \frac{\Delta}{\ell}\right)=1$, then $\ell$ is called an \textit{Elkies prime}. In this case, the characteristic equation factors as $\chi(\phi_p)=(\phi_p-\lambda)(\phi_p-\mu) \mod \ell$, so when acting on $E[\ell]$ the map $\phi_p$ has two eigenvalues $\lambda, \mu \in \F_\ell^*$ with corresponding eigenpoints $P, Q$. Since $\lambda\mu=p$ and $\lambda+\mu=t$, it obviously suffices to determine one of them. So we have to solve the discrete logarithm problem
\[ 
\lambda P=\phi_p(P)=(P_x^p, P_y^p),
 \]
which only requires working in extensions of degree $O(\ell)$.
\item If  $\left( \frac{\Delta}{\ell}\right)=-1$, then $\ell$ is called an \textit{Atkin prime}. In this case the eigenvalues of $\phi_p$ are in $\F_{\ell^2}\backslash\F_\ell$ and there is no eigenpoint $P \in E[\ell]$. There is a generic method for computing the value of $t \mod \ell$ for Atkin primes, which is of equal run-time as the one available for Elkies primes. However, it does not yield the exact value of $t \mod \ell$ but only a set of candidates and is thus only efficient provided the cardinality of this set is small.
\end{enumerate}
\section{Polynomially cyclic algebras}\label{sec:poly_cyc_alg}
In this section we recall facts on polynomially cyclic algebras, which were first described in \cite{MiVu}. This general framework will be used in the next section in order to elaborate a new approach to the Atkin case in Schoof's algorithm.\\

\begin{defi}\cite[p.~6]{MiVu}
Let $\mathbb{K}$ be a finite field and $f(X) \in \mathbb{K}[X]$ a polynomial with $\deg2(f)=n$. Then we call the $\mathbb{K}$-algebra $\mathbf{A}=\mathbb{K}[X]/(f(X))$ a \textit{polynomially cyclic algebra} with \textit{cyclicity polynomial} $C(X) \in \mathbb{K}[X]$ and we call $f(X)$ a \textit{cyclic polynomial} if the following conditions are satisfied:
\begin{enumerate}
\item $f(C(X)) \equiv 0 \textup{ mod } f(X)$,
\item $C^{(n)}(X)-X \equiv 0 \textup{ mod } f(X)$ and $\gcd(C^{(m)}(X)-X, f(X))=1$ for $m<n$.
\end{enumerate}
Here $C^{(m)}(X)=\underbrace{C \circ C \dots \circ C}_{m\textup{ times}}(X)$.
\end{defi}
\begin{bem}
\begin{enumerate}
\item Every irreducible polynomial $f \in \mathbb{K}[X]$ is cyclic as well. Its cyclicity polynomial is $C(X)=X^q$, where $\mathbb{K}=\F_q$.
\item If a polynomial $f(X) \in \mathbb{K}[X]$ is cyclic with cyclicity polynomial $C(X)$, the same holds true in all extensions $\mathbb{L}$ of $\mathbb{K}$.
\end{enumerate}
\end{bem}

The following theorem provides several ways to describe cyclic polynomials.

\begin{satz}\cite[p.~6]{MiVu}
The following conditions are equivalent:
\begin{enumerate}
\item $f$ is cyclic.
\item There exists a polynomial $C(X) \in \mathbb{K}[X]$ which cyclically permutes the roots of $f$, i.~e., for every root $\alpha \in \overline{\mathbb{K}}$ of $f$ the equality $f(C(\alpha))=0$ holds and the elements $C^{(i)}(\alpha),$ $i=1, \ldots, n$, are pairwise distinct.
\item In the factorisation $f=\prod_{k=1}^d h_k$ over $\mathbb{K}[X]$ all factors $h_k$ have the same degree and are pairwise distinct.
\end{enumerate}
\end{satz}

Next we describe some of the properties of polynomially cyclic algebras.

\begin{satz}\cite[p.~9]{MiVu}\label{satz:automorphismen}
Let $\mathbf{A}=\mathbb{K}[X]/(f(X))$, $\deg2(f(X))=n$, be a polynomially cyclic algebra with cyclicity polynomial $C(X)$, let $\alpha:=X \mod f(X)$. Then the following statements hold:
\begin{enumerate}
\item The cyclicity polynomial $C(X)$ induces an automorphism $\nu$ of $\mathbb{K}$-algebras of order $n$ in virtue of
$$\nu: \mathbf{A} \rightarrow \mathbf{A},\quad \alpha \mapsto C(\alpha).$$

We write $\gal(\mathbf{A}/\mathbb{K}):=\langle\nu \rangle$, so $\gal(\mathbf{A}/\mathbb{K})$ is the automorphism group of $\mathbf{A}$, generated by $\nu$, which we also call the galois group of $\mathbf{A}/\mathbb{K}$.

\item We have $\mathbf{A}^{\gal(\mathbf{A}/\mathbb{K})}=\textup{Fix}(\nu)=\mathbb{K}$, so $\nu(x)=x \Rightarrow x \in \mathbb{K}$ holds.
\end{enumerate}
\end{satz}

\begin{satz}\cite[p.~11]{MiVu}\label{satz:unteralgebren}
Let $\mathbf{A}=\mathbb{K}[X]/(f(X))$ be a polynomially cyclic algebra. Then the following statements hold:
\begin{enumerate}
\item Let $\tilde{\mathbb{K}}/\mathbb{K}$ be a field extension. Then $\tilde{\mathbf{A}}:=\mathbf{A} \otimes_\mathbb{K} \tilde{\mathbb{K}}$ is a polynomially cyclic $\tilde{\mathbb{K}}$-algebra and there is a canonical isomorphism $\gal(\tilde{\mathbf{A}}/\tilde{\mathbb{K}}) \cong \gal(\mathbf{A}/\mathbb{K})$.

\item Let $H \subset \gal(\mathbf{A}/\mathbb{K})$ be a subgroup. Then the subalgebra of all elements invariant under $H$,
$$\mathbf{A}^H:=\{a \in \mathbf{A}: h(a)=a \ \forall h \in H\},$$
is polynomially cyclic. Conversely, if $\mathbf{B} \subset \mathbf{A}$ is a polynomially cyclic algebra there exists $H \subset \gal(\mathbf{A}/\mathbb{K})$ such that $\mathbf{B}=\mathbf{A}^H$.\\
The dimension of $\mathbf{A}^H$ equals the index $[\gal(\mathbf{A}/\mathbb{K}):H]$. In addition, there is a canonical isomorphism $\gal(\mathbf{A}^H/\mathbb{K}) \cong \gal(\mathbf{A}/\mathbb{K})/H$.\\
\end{enumerate}
\end{satz}

Starting from our definitions and using the properties of polynomially cyclic algebras we have mentioned we now proceed to define \textit{Lagrange resolvents} in these. Again, we closely follow \cite[p.~12]{MiVu}. Let $\mathbf{A}=\mathbb{K}[X]/(f(X))$, $\deg2(f(X))=n$, again be a polynomially cyclic algebra with $\gal(\mathbf{A}/\mathbb{K})=\langle \nu \rangle$. Let $\rho_n \in \overline{\mathbb{K}}$ be a primitive $n$-th root of unity with minimal polynomial  $K_n(X) \in \mathbb{K}[X]$. We now define
$$\mathbf{A}_{\rho_n} = \mathbf{A}[T]/(K_n(T)).$$
Since $\mathbf{A}_{\rho_n} =  \mathbf{A} \otimes_\mathbb{K} \mathbb{K}[\rho_n]$ holds, theorem \ref{satz:unteralgebren} implies $\mathbf{A}_{\rho_n}$ is a polynomially cyclic algebra over $\mathbb{K}[\rho_n]$ and using the canonical isomorphism we can identify the groups $\gal(\mathbf{A}_{\rho_n}/ \mathbb{K}[\rho_n])$ and  $\gal(\mathbf{A}/ \mathbb{K})$.
Now let
$\chi: \gal(\mathbf{A}/ \mathbb{K}) \rightarrow \mathbb{K}[\rho_n]$ 
be a multiplicative character. For $\alpha \in \mathbf{A}$ we define the Lagrange resolvent $(\chi, \alpha)$ as
$$(\chi, \alpha)= \sum_{\sigma \in \gal(\mathbf{A}/\mathbb{K})} \chi(\sigma)\sigma(\alpha)
= \sum_{i=1}^n \chi(\nu^i)\nu^i(\alpha)=  \sum_{i=1}^n \chi(\nu)^i\nu^i(\alpha) \in \mathbf{A}_{\rho_n}.$$

The following theorem establishes some properties of the Lagrange resolvent.


\begin{satz}\label{satz:lagrange-eigenschaften} \cite[pp.~12-13]{MiVu}
Let $\chi, \chi': \gal(\mathbf{A}/ \mathbb{K}) \rightarrow \mathbb{K}[\rho_n]$ be two characters. Then the following statements hold:

\begin{enumerate}
\item For $\sigma \in \gal(\mathbf{A}/ \mathbb{K})$ we obtain
$$\sigma(\chi, \alpha)= \chi^{-1}(\sigma)(\chi, \alpha).$$

\item If $r=\ord(\chi)$, then
$$(\chi, \alpha)^r \in \mathbb{K}[\rho_n].$$

\item Likewise, we obtain
$$\frac{(\chi, \alpha) \cdot (\chi', \alpha)}{(\chi \cdot \chi', \alpha)} \in \mathbb{K}[\rho_n],\quad \text{provided}\quad (\chi \cdot \chi', \alpha) \in \mathbf{A}^* 
\quad \text{holds}.$$

\item For $\mathbb{K}=\F_q$ and $\varphi_q \in \gal(\mathbf{A}/ \mathbb{K}),\ x \mapsto x^q$ we have
$$(\chi, \alpha)^q = \chi^{-q}(\varphi_q)(\chi^q, \alpha).$$
\end{enumerate}
\end{satz}

The Lagrange resolvents can now be used to solve the following general problem, which we will afterwards consider in a special case:\\
Let $\mathbf{A}_i=\mathbb{K}[X]/(f_i(X))$, $\deg2(f_i(X))=n$, with cyclicity polynomials $C_i(X)$, $i=1, 2$, be two isomorphic polynomially cyclic algebras. Denoting by $\nu_i$ the automorphism induced by $C_i(X)$ and setting $G_i:=\gal(\mathbf{A}_i/\mathbb{K})=\langle \nu_i\rangle$ as well as $\alpha_i:=X \mod f_i(X)$, we wish to determine, more precisely, an isomorphism
\begin{equation}\label{eq:iso_allgemein}
\varphi: \mathbf{A}_1 \rightarrow \mathbf{A}_2,\quad \alpha_1 \mapsto \sum_{i=1}^n a_i\nu_2^i(\alpha_2)\quad \text{with}\quad a_i \in \mathbb{K}
\end{equation}
such that $\varphi \circ \nu_1 = \nu_2 \circ \varphi$ holds.
\begin{bem}\label{bem:normale_basis}
We can only obtain an isomorphism in this shape if $\alpha_2$ forms a normal basis of $\mathbf{A}_2$ together with its conjugates. It turns out, however, that for the practical application we consider this seems to be always the case, which yields the existence of isomorphism \eqref{eq:Iso} below.
\end{bem}

In order to determine isomorphism \eqref{eq:iso_allgemein} we can avail ourselves of the following statement:

\begin{satz}\cite[p.~14]{MiVu}\label{satz:iso_koeffizienten_allgemein}
Let, as just mentioned, $\mathbf{A}_i$ be isomorphic polynomially cyclic $\mathbb{K}$-algebras and $\varphi$ be the isomorphism 
$\varphi: \mathbf{A}_1 \tilde{\rightarrow} \mathbf{A}_2$. 
Let further $\rho_n \in \overline{\mathbb{K}}$ be a primitive $n$-th root of unity and $\mathbf{A}_{i, \rho_n}$ be defined as above. Let $\chi_i: G_i \rightarrow \langle \rho_n \rangle,$ $i=1, 2,$ be characters with $\chi_1(\nu_1)=\chi_2(\nu_2)$. Then there exists $\beta(\chi_2) \in \mathbb{K}[\rho_n]$, such that
$$\varphi((\chi_1, \alpha_1))=(\chi_2, \alpha_2)\cdot \beta(\chi_2).$$
More precisely,
$$\beta(\chi_2)=\sum_{i=1}^n a_i\chi_2^{-1}(\nu_2^i)$$
holds.
\end{satz}
By means of this theorem the coefficients  $a_i$ in \eqref{eq:iso_allgemein} can be determined in the following way \cite[p.~14]{MiVu}:\\
For $j=1, \ldots, n$ let $\chi_{2,j}$ be a character with $\chi_{2,j}(\nu_2)=\rho_n^j$ and assume the value $\beta(\chi_{2,j})$ is known. This implies
$$\beta(\chi_{2,j})=\sum_{i=1}^n a_i\chi_{2,j}^{-1}(\nu_2^i)
=\sum_{i=1}^n a_i\rho_n^{-ij},\textup{ }j=1, \ldots, n.$$
Hence,
$$M \cdot \vec{a} = \vec{\beta},$$
where $M=(\rho_n^{-ji})_{j,i=1}^{n}$, $\vec{a}=(a_i)_{i=1}^n$ and $\vec{\beta}=(\beta(\chi_{2,j}))_{j=1}^n$. One thus obtains a linear system of equations for the coefficients $a_i$, which can be used to determine them because of the following

\begin{prop}\label{prop:matrix_regulaer}
The matrix $M$ is regular if $\gcd(n,char(\mathbb{K}))=1$ holds. More precisely, 
$$(\det(M))^2=(-1)^{n \cdot (n+1)/2 +1} \cdot n^n$$
holds.
\begin{proof}
It is obvious that the columns of $M$ may be rearranged to form a Vandermonde matrix of the form $M'=(\rho_n^{ij})_{i, j=0}^{n-1}$ where $M'_{ij}=(\rho_n^i)^j$ and $\det(M')^2=\det(M)^2$ holds. The general formula for the determinant of Vandermonde matrices implies
$$\det(M')=\prod_{1\leq i <j \leq n} (\rho_n^j-\rho_n^i),$$
hence
$$
\begin{aligned}
(-1)^{n \cdot (n-1)/2}\det(M')^2=& \prod_{1\leq i \neq j \leq n} (\rho_n^j-\rho_n^i)
=\prod_{j=1}^n \prod_{\substack{i=1 \\ i \neq j}}^n (\rho_n^j-\rho_n^i)\\
=& \prod_{j=1}^n \rho_n^j \prod_{\substack{i=1 \\ i \neq j}}^n (1-\rho_n^{i-j})
=\prod_{j=1}^n \rho_n^j \prod_{i=1}^{n-1} (1-\rho_n^i).
\end{aligned}
$$
Since $X^{n}-1=(X-1)(X^{n-1} + \cdots +X+1)$,  $X^{n-1} + \cdots +X+1=\prod_{i=1}^{n-1} (X-\rho_n^i)$ holds, one obtains $\prod_{i=1}^{n-1} (1-\rho_n^i)=n$. Hence,
$$(-1)^{n \cdot (n-1)/2} \det(M')^2=\prod_{j=1}^n (\rho_n^j \cdot n)=n^n \rho_n^{\sum_{j=1}^n j}
=n^n \rho_n^{n\cdot (n+1)/2}=(-1)^{n+1} \cdot n^n,
$$
so $\det(M')^2=(-1)^{n \cdot (n+1)/2 +1} \cdot n^n$.
\end{proof}
\end{prop}
\section{Application}\label{sec:anwendung}
We first recall the ray-polynomial from \cite{MiVu}. 

\begin{defi}\label{def:defi_strahlenpolynom}
    Let $E$ be an elliptic curve over $\F_p$, $\ell$ a prime and $P$ a point in  $E[\ell], P \neq \OO$. Then the \textit{ray-polynomial} corresponding to $P$ is defined by
    $$E_P(X)=\prod_{a=1}^{(\ell-1)/2} (X - (aP)_x) \in \overline{\F}_p[X].$$
\end{defi}
We remark that the ray-polynomial depends only on the subspace of $E[\ell]$ spanned by $P$. It has the following properties:
\begin{lemma}
\begin{enumerate}
\item $E_P(X)$ is a cyclic polynomial, whose cyclicity polynomial $G_c(X)$ can be easily computed from the well-known division polynomials of $E$ \cite[p.~8]{MiVu}.
\item The field of definition of $E_P(X)$ is $\F_{p^r}$, where $r$ is the degree of an irreducible factor of the modular polynomial $\Phi_\ell(X, j(E)) \in \F_p[X]$ (for a definition cf.~\cite{Cox}). This follows from \cite[Theorems 6.1, 6.2]{Schoof} and is shown in \cite[p.~3]{MiVu}. In the Elkies case, for an appropriate $P$ the ray-polynomial coincides with the Elkies factor $f_{\ell, \lambda}$ (cf.~\cite{Schoof, Morain}) and $r$ is thus $1$.
\end{enumerate}
\label{lem:strahlen_poly_eig}
\end{lemma}

Our overall strategy is the same as in Schoof's algorithm. We wish to determine the value $t \mod \ell$ by considering the equation $\chi(\phi_p) \mod \ell$. Plugging in an $\ell$-torsion point $P$ and restricting to $x$-coordinates we obtain
\begin{equation}\label{eq:char_gleichung_x}
(\phi_p^2(P)+pP)_x=G_t(\varphi_p(P_x)),
\end{equation}
where $G_t$ is as in equation \eqref{eq:multi_poly}. Setting $\An=\F_{p^r}$ as in lemma \ref{lem:strahlen_poly_eig} and $\Bn=\An[T]/(E_P(T))$ all computations can be performed in $\Bn$. Lemma \ref{lem:strahlen_poly_eig} again implies that $\Bn$ is a polynomially cyclic algebra. We denote by $\nu$ the generator of the galois group $\gal(\Bn/\An)$ induced by the cyclicity polynomial $G_c$. As in Schoof's algorithm, calculating the action of $\phi_p$ dominates the run-time. Since $\deg2(E_P)=(\ell-1)/2$, the overall complexity comprises $O(r\ell \log p)$ operations in $\F_p$.\\
Using the results presented in section \ref{sec:poly_cyc_alg} we want to describe an approach allowing to decrease the run-time. Our main idea, which was first sketched in \cite{MiVu} and worked out in detail in \cite{ma}, is to construct a polynomially cyclic algebra $\Cn$ which is isomorphic to the algebra $\Bn$ defined by means of $E_P$ and which allows for an efficient computation of the Frobenius homomorphism. After that, we wish to solve the resulting discrete logarithm problem in this algebra. Our approach is applicable for Atkin primes and contrasts with the various improvements available for the Elkies case, which essentially rely on transferring computations into smaller extensions of $\F_p$, a strategy which is impossible for Atkin primes. Obviously, our approach requires that the isomorphism between the two algebras be explicitly computed.\\

\subsection{The isomorphic algebra}
We now proceed to define the algebra $\Cn$.

\begin{prop}\label{prop:defi_C}
    Let $G=\left\{b \in \F_\ell^*: \left(\frac{b}{\ell} \right)=1\right\}$, $\zeta_\ell \in \overline{\F}_p$ be a primitive $\ell$-th root of unity and $K(U)=\prod_{b \in G} (U-\zeta_\ell^b)$. Then $\mathbf{C}=\mathbf{A}[U]/(K(U))$ is a polynomially cyclic algebra with $\gal(\mathbf{C}/\mathbf{A})=\langle \sigma: \zeta_\ell \mapsto \zeta_\ell^{c^2} \rangle$, where $\langle c \rangle=\F_\ell^*$ holds.
    \begin{proof} We first prove that $K(U)$ indeed lies in $\mathbf{A}[U]$. From \cite[Theorem 6.1, 6.2]{Schoof} we deduce that $\phi_p^r$ acts on $E[\ell]$ as a scalar matrix $\left( \begin{smallmatrix} a & 0 \\ 0 & a
\end{smallmatrix} \right)$. Hence, $\phi_p^r$ exhibits the double eigenvalue $a$. Denoting $\lambda, \mu$ the eigenvalues of $\phi_p$ this implies $\lambda^r=a=\mu^r$. Since $\lambda\mu=p \mod \ell$ this implies $p^r=a^2 \mod \ell$ and thus $\left(\frac{p^r}{\ell} \right)=1$. This means the Frobenius homomorphism $\varphi_{p^r}: \mathbf{A} \rightarrow \mathbf{A},\ x \mapsto x^{p^r}$ maps roots of $K(U)$ to other roots of this polynomial, which thus lies in $\mathbf{A}[U]$.\\
Now let $c$ be a generator of $\F_\ell^*$. Then obviously $c^2$ generates $G$, which implies the polynomial $C(X)=X^{c^2}$ permutes the roots of $K$. So $K$ is cyclic.
    \end{proof}
\end{prop}

\begin{lemma}
    Using the above notations $\mathbf{B} \cong \mathbf{C}$ as algebras over $\mathbf{A}$.
    \begin{proof}
	First, we know
        $$\deg2(K(U))=\#(\F_\ell^*)^2=\frac{\ell-1}{2}=\deg2(E_P(T)).$$
	Furthermore, according to the above considerations $\mathbf{B}$ as well as $\mathbf{C}$ are polynomially cyclic algebras and hence the polynomials $K(U)$ as well as $E_P(T)$ decompose into irreducible factors of equal degree. It now remains to show that the degree of the factors in the factorisation of the two polynomials coincides, which implies that the factors occurring in the decomposition of $\Bn$ and $\Cn$ into a product of fields are isomorphic. This follows since we work over $\F_p$ and it is thus sufficient to show that these fields have the same degree over $\F_p$.
    Since $\mathbf{A}=\F_{p^r}$, $[\mathbf{A}:\F_p]=r$ holds and thus for $\mathbb{L} \supset \mathbf{A}$ we have
        $$[\mathbb{L}:\mathbf{A}]=\min\{m: \varphi_p^{rm}(z)=z \ \forall z \in \mathbb{L} \}.$$
        We first consider $\mathbf{B}$. As mentioned in proposition \ref{prop:defi_C} there is $a \in \F_\ell^*$ such that $\phi_p^r(P)=aP$ for all $P \in E[\ell]$. Now let  $\Ln_1$ be a factor of $\mathbf{B}$. Writing $\theta:=T+(E_P(T))$, we obtain
        $$
        \begin{aligned}
        \left[\Ln_1:\mathbf{A}\right]
        &=\min\{m: \varphi_p^{rm}(z)=z \ \forall z \in \Ln_1\}
        =\min\{m: \varphi_p^{rm}(\theta)=\theta\}\\
        &=\min\{m: (a^mP)_x=\theta\}
        =\min\{m: a^mP=\pm P\}\\
        &=\ord_S(a),
        \end{aligned}
        $$
        
        where $S=\F_\ell^*/\{\pm 1\}$.\\
        
        Now let $\Ln_2$ be a factor of $\mathbf{C}$. Anew we use the fact that $a^2\equiv p^r \mod \ell$  and hence $p^r \in G$ holds. It follows
        
        $$
        \begin{aligned}
        \left[\Ln_2:\mathbf{A}\right]
        &=\min\{m: \varphi_p^{rm}(z)=z \ \forall z \in \Ln_2\}
        =\min\{m: \varphi_p^{rm}(\zeta_\ell)=\zeta_\ell\}\\
        &=\min\{m:\zeta_\ell^{p^{rm}}=\zeta_\ell\}
        =\min\{ m: p^{rm} \equiv 1 \mod \ell \}\\
        &=\ord_G(p^r).
        \end{aligned}
        $$
        
        As the groups $S$ and $G$ are isomorphic by virtue of $c \mapsto c^2$, we glean
        $$[\Ln_1:\mathbf{A}]=\ord_S(a)=\ord_G(a^2)=\ord_G(p^r)=[\Ln_2:\mathbf{A}].$$
    \end{proof}
\end{lemma}

So there exists an isomorphism
\begin{equation}\label{eq:Iso}
\alpha: \mathbf{B} \rightarrow \mathbf{C},\quad \theta \mapsto \sum_{i=1}^{(\ell-1)/2}b_i\zeta_\ell^{c^{2i}}\quad \text{with}\quad b_i \in \mathbf{A}.
\end{equation}
Furthermore, we require $\alpha \circ \nu = \sigma\circ \alpha$. The isomorphism should thus commute with the automorphism of both algebras such that the prerequisites of theorem \ref{satz:iso_koeffizienten_allgemein} are satisfied.

\begin{bem}
As already stated in remark \ref{bem:normale_basis} we only obtain an isomorphism of shape \eqref{eq:Iso} if $\zeta_\ell^{c^{2i}}, i=1, \ldots, \frac{\ell-1}{2}$, form a basis of the algebra $\mathbf{C}$. This can easily be checked during actual computations and has in practice always been the case.
\end{bem}

\subsection{Construction of the isomorphism}\label{sec:iso_rechnung}
In order to determine the coefficients of the isomorphism \eqref{eq:Iso} we follow the explanations after theorem \ref{satz:iso_koeffizienten_allgemein}.\\

Let $\rho$ be an $(\ell-1)/2$-th root of unity, $q\mid (\ell-1)/2$, $\rho_q=\rho^{(\ell-1)/(2q)}$ and  
\begin{equation}\label{eq:def_chi_q}
\chi_q: \gal(\mathbf{B}/\mathbf{A}) \rightarrow \mathbf{A}[\rho_q],\text{ } \nu \mapsto \rho_q
\end{equation} a character of order $q$. We identify $\chi_q$ with the character
$$\chi_{q, 2}: \gal(\mathbf{C}/\mathbf{A}) \rightarrow \mathbf{A}[\rho_q],\text{ } \sigma \mapsto \rho_q.$$   Then the following holds:
\begin{lemma}
    Let $$b^{(q)}_i=\sum_{k=1}^{(\ell-1)/(2q)} b_{kq+i}, 1\leq i \leq q, \text{ and } \theta^{(q)}=\sum_{j=1}^{(\ell-1)/(2q)}\nu^{{jq}}(\theta) \text{ and }
    \zeta_\ell^{(q)}=\sum_{j=1}^{(\ell-1)/(2q)}\sigma^{jq}(\zeta_\ell).$$
    Let further $$\tau_e(\chi_q)=\sum_{a=1}^{q}\rho_q^a \nu^{a}(\theta^{(q)}) 
    \text{ and } \tau(\chi_q)=\sum_{j=1}^{q}\rho_q^{j}\sigma^{j}(\zeta_\ell^{(q)})
    \text{ as well as } \beta(\chi_q)=\sum_{i=1}^{q}\rho_q^{-i}b^{(q)}_i.$$
    Then 
    $$\alpha(\tau_e(\chi_q))=\beta(\chi_q) \cdot \tau(\chi_q).$$
    
    \begin{proof}
		Using the isomorphism $\alpha$ and the property $\alpha \circ \nu = \sigma \circ \alpha$ we compute
        
        \begin{align*}
        \alpha((\chi_q, \theta))&=\sum_{a=1}^{(\ell-1)/2}\rho_q^a \alpha(\nu^{a}(\theta))
        =\sum_{a=1}^{(\ell-1)/2}\rho_q^a\sigma^a\left(\sum_{i=1}^{(\ell-1)/2}b_i\zeta_\ell^{c^{2i}}\right)\\
        &=\sum_{a=1}^{(\ell-1)/2}\rho_q^a \sum_{i=1}^{(\ell-1)/2}b_i\sigma^a(\zeta_\ell^{c^{2i}})
        =\sum_{i=1}^{(\ell-1)/2}b_i \rho_q^{-i} \sum_{a=1}^{(\ell-1)/2}\rho_q^{a+i}\sigma^{a+i}(\zeta_\ell)\\
        &=\sum_{i=1}^{(\ell-1)/2}b_i \rho_q^{-i} \sum_{j=1}^{(\ell-1)/2}\rho_q^{j}\sigma^{j}(\zeta_\ell)
        =\sum_{i=1}^{q}\rho_q^{-i} \sum_{k=1}^{(\ell-1)/(2q)} b_{kq+i} \sum_{j=1}^{(\ell-1)/2}\rho_q^{j}\sigma^{j}(\zeta_\ell)\\
        &=\beta(\chi_q)\cdot (\chi_q, \zeta_\ell).
        \end{align*}
        
		Using the above definitions we first glean
        $\beta(\chi_q)=\sum_{i=1}^{q}\rho_q^{-i}b^{(q)}_i$.
        Further, we obtain 
        $$(\chi_q, \theta)=\sum_{a=1}^{(\ell-1)/2}\rho_q^a\nu^{a}(\theta)
        =\sum_{a=1}^{q} \sum_{j=1}^{(\ell-1)/(2q)} \rho_q^{jq+a} \nu^{{jq+a}}(\theta)
        =\sum_{a=1}^{q}\rho_q^a \nu^{a}(\theta^{(q)}).$$ 
        Similarly, we get $(\chi_q, \zeta_\ell)=\sum_{j=1}^{q}\rho_q^{j}\sigma^{j}(\zeta_\ell^{(q)})$.    
        
    \end{proof}
\end{lemma}

\begin{bem}
The quantities $\tau(\chi_q)$ are essentially 
the well-known cyclotomic Gau{\ss} sums. Since they are formed in natural analogy to these, the values $\tau_e(\chi_q)$ were named \textit{elliptic Gau{\ss} sums} in \cite{MiVu}.
\end{bem}

We first concern ourselves with the computation of $\beta(\chi_q)$. The fact that $\ord(\chi_q)=q$ holds and theorem \ref{satz:lagrange-eigenschaften} imply that $\tau(\chi_q)^q$ as well as $\tau_e(\chi_q)^q$ lie in $\mathbf{A}[\rho_q]$. From the last lemma we deduce 
$$\tau_e(\chi_q)^q=\alpha(\tau_e(\chi_q)^q)=\alpha(\tau_e(\chi_q))^q=\beta(\chi_q)^q\cdot \tau(\chi_q)^q.$$ Hence, we obtain
\begin{equation}\label{eq:berechne_chi_q}
\beta(\chi_q)^q=\frac{\tau_e(\chi_q)^q}{\tau(\chi_q)^q}.
\end{equation}

The cost for computing $\beta(\chi_q)$, which lies in $\mathbf{A}[\rho_q]$ by definition, thus consists in calculating the $q$-th powers of the two Lagrange resolvents and in extracting a $q$-th root in $\mathbf{A}[\rho_q]$.\\

Having determined $\beta(\chi_q)$ for a character of order $q$, we can employ another one of the properties from theorem \ref{satz:lagrange-eigenschaften} to compute $\beta(\chi_q^i),$ $i=2, \ldots, q$. Namely, property 3 implies
$$z_{e,1}:=\frac{\tau_e(\chi_q)\cdot \tau_e(\chi_q)}{\tau_e(\chi_q^2)} \in \mathbf{A}[\rho_q] \quad
 \text{and}\quad z_1:=\frac{\tau(\chi_q)^2}{\tau(\chi_q^2)} \in \mathbf{A}[\rho_q].$$
Likewise 
$$z_{e,i}:=\frac{\tau_e(\chi_q^i)\cdot \tau_e(\chi_q)}{\tau_e(\chi_q^{i+1})} \in \mathbf{A}[\rho_q] \quad
\text{and}\quad z_i:=\frac{\tau(\chi_q^i) \cdot \tau(\chi_q)}{\tau(\chi_q^{i+1})} \in \mathbf{A}[\rho_q]$$
holds. This yields 
\begin{equation}\label{eq:berechne_pot_von_beta_chi_q}
\beta(\chi_q^2)=\frac{\alpha(\tau_e(\chi_q^2))}{\tau(\chi_q^2)}=\frac{\alpha(\tau_e(\chi_q)^2)}
{\tau(\chi_q)^2}\cdot \frac{z_{1}}{z_{e,1}}=\beta(\chi_q)^2\cdot \frac{z_1}{z_{e,1}}\quad
\text{and} \quad \beta(\chi_q^{i+1})=\beta(\chi_q^i)\beta(\chi_q)\cdot \frac{z_i}{z_{e,i}}.
\end{equation}

Hence, the values $\beta(\chi_q^i),$ $i=2, \ldots, q,$ can be determined without extracting a root again.\\

Now assume that for all maximal prime divisors $q \mmid (\ell-1)/2=:n$ the values $\beta(\chi_q^{i_q}),$ $i_q=1, \ldots, q,$ have been determined and let $\chi_n$ be a character of order $n$, so
$\chi_n=\prod_{q \mmid n} \chi_q^{e_q}$. Then we obtain
$$\beta(\chi_n)=\frac{\alpha(\tau_e(\chi_n))}{\tau(\chi_n)}=\frac{\alpha\left(\tau_e\left(
\prod_{q \mmid n} \chi_q^{e_q}\right)\right)}{\tau\left(\prod_{q \mmid n} \chi_q^{e_q}\right)}.$$
Here,
$$\alpha\left(\tau_e\left(\prod_{q \mmid n} \chi_q^{e_q}\right)\right)=
\prod_{q \mmid n} \alpha(\tau_e(\chi_q^{e_q})) \cdot \frac{\alpha\left(\tau_e\left(\prod_{q \mmid n} \chi_q^{e_q}\right)\right)}{\prod_{q \mmid n} \alpha(\tau_e(\chi_q^{e_q}))}=:\prod_{q \mmid n} \alpha(\tau_e(\chi_q^{e_q})) \cdot z_{e,n}$$
and in a similar vein $\tau\left(\prod_{q \mmid n} \chi_q^{e_q}\right)=\prod_{q \mmid n} \tau(\chi_q^{e_q}) \cdot z_{n}$, where $z_{e,n}, z_n \in \mathbf{A}[\rho]$ holds.
Thus, we derive the equation
\begin{equation}\label{eq:beta_allgemein_berechnen}
\beta(\chi_n)=\frac{\prod_{q \mmid n} \alpha(\tau_e(\chi_q^{e_q})) \cdot z_{e,n}}
{\prod_{q \mmid n} \tau(\chi_q^{e_q}) \cdot z_{n}}=\frac{z_{e,n}}{z_n}\cdot \prod_{q\mmid n} \beta(\chi_q^{e_q}).
\end{equation}

Choosing $\chi_q$ for all $q \mmid n$  as in \eqref{eq:def_chi_q} as a primitive character allows us  to compute $\beta(\chi)$ for all characters $\chi: \gal(\mathbf{B}/\mathbf{A}) \rightarrow \mathbf{A}[\rho]$ by evaluating this formula for $e_q=1, \ldots, q$.\\

Having computed these values, we find ourselves in the situation described by theorem \ref{satz:iso_koeffizienten_allgemein} and are provided with a linear system of the form
\begin{equation}\label{eq:iso_bestimmung_lgs}
M\cdot \vec{b}=\vec{\beta},
\end{equation}
where $M=(\rho_n^{-ji})_{j,i=1}^{n}$, $\vec{b}=(b_i)_{i=1}^n$ and $\vec{\beta}=(\beta(\chi))$. Since $\gcd(n, p)=1$, the matrix $M$ is regular according to proposition \ref{prop:matrix_regulaer}. Hence, we can finally determine the coefficients $b_i$ of the isomorphism $\alpha$.

\subsubsection{Improvement}\label{sec:alpha_best_besser}
In order to avoid computing $\beta(\chi)$ for all the characters $\chi$ of order $n=\frac{\ell-1}{2}$, which produces major costs, in this section we present an alternative approach for determining $\alpha$.\\

First, using the definition of $\alpha$ we derive
\[ 
\alpha(\theta)=\sum_{i=1}^{(\ell-1)/2} b_i \sigma^i(\zeta_\ell)\quad \text{and}\quad \alpha(\nu^k(\theta))=\sum_{i=1}^{(\ell-1)/2} b_i \sigma^{i+k}(\zeta_\ell).
 \]
By means of these identities we calculate
\begin{align*}
\alpha(\theta^{(q)})=&\alpha\left(\sum_{j=1}^{(\ell-1)/(2q)} \nu^{jq}(\theta)\right)
=\sum_{j=1}^{(\ell-1)/(2q)} \sum_{i=1}^{(\ell-1)/2} b_i \sigma^{i+jq}(\zeta_\ell)\\ 
=& \sum_{i=1}^{(\ell-1)/2} b_i \sigma^i \left( \underbrace{\sum_{j=1}^{(\ell-1)/(2q)} \sigma^{jq}(\zeta_\ell) }_{\zeta_\ell^{(q)}} \right )
=\sum_{i=1}^{q} \underbrace{\sum_{k=1}^{(\ell-1)/(2q)} b_{kq+i}}_{b_i^{(q)}} \underbrace{\sigma^{kq+i}(\zeta_\ell^{(q)})}_{=\sigma^i(\zeta_\ell^{(q)})}
=\sum_{i=1}^{q} b_i^{(q)}\sigma^i(\zeta_\ell^{(q)}).
\end{align*}

Hence, we require exactly the values $b_i^{(q)}$ to specify the isomorphism
\[ 
\alpha_q: \An[\theta^{(q)}] \rightarrow \An[\zeta_\ell^{(q)}]
 \]
arising by restriction of $\alpha$ to these sub-algebras.\\
To determine the $b_i^{(q)}$ one has to compute the values $\beta(\chi_q)$ by extracting one $q$-th root as in section \ref{sec:iso_rechnung}. Subsequently, one directly proceeds to solve a linear system of equations instead of determining $\beta(\chi)$ for general characters of order $n$. Our new approach consists in computing the isomorphism $\alpha_q$ for $q \mmid \frac{\ell-1}{2}$ and in inductively constructing the isomorphism $\alpha$ from these intermediate data.\\

Let $q_1, q_2 \mid \frac{\ell-1}{2}$ and $(q_1, q_2)=1$ and assume the isomorphisms $\alpha_{q_1}, \alpha_{q_2}$ have been determined. We present a procedure to  compute $\alpha_{q_1q_2}$. First, set $\Bn_q:=\An[\theta^{(q)}], \Cn_q:=\An[\zeta_\ell^{(q)}]$ and consider the following diagram:

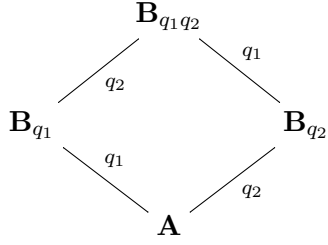
\begin{displaymath}
    \xymatrix{ &\Bn_{q_1q_2} \ar@{-}[ld]^{q_2} \ar@{-}[rd]^{q_1}& \\
              \Bn_{q_1} \ar@{-}[rd]^{q_1}& & \Bn_{q_2} \ar@{-}[ld]^{q_2} \\
              &\An & }
\end{displaymath}

Using the general theory from \cite{MiVu} we obtain $\gal(\Bn/\Bn_q)=\langle \nu^q \rangle$, where $\gal(\Bn/\An)=\langle \nu \rangle$ holds. The polynomial
\[
M_1(X)=\prod_{i=1}^{q_2} X-\nu^{q_1i}(\theta^{(q_1q_2)})
\]
vanishes at $\theta^{(q_1q_2)}$, and since its roots are obviously permuted by $\nu^{q_1}$ it lies in $\Bn_{q_1}[X]$. Further,
\[
M_2(X)=\prod_{i=1}^{q_2} X-\nu^{i}(\theta^{(q_2)}) 
\]
has $\theta^{(q_2)}$ as a root. Considering the action of $\nu^{q_1}$ on the roots of $M_2(X)$ we observe
\[ 
\nu^{jq_1}(\nu^i(\theta^{(q_2)}))=\nu^{jq_1+i}(\theta^{(q_2)})=\nu^{j(i)}(\theta^{(q_2)}),
 \]
where $j(i)\equiv jq_1+i \mod q_2$, as $\nu^{q_2}(\theta^{(q_2)})=\theta^{(q_2)}$ holds. Since $(q_1, q_2)=1$, one deduces $jq_1 \nequiv 0 \mod q_2$ for $j<q_2$. Thus, $\nu^{q_1}$ permutes the roots of $M_2(X)$, which thus lies in $\Bn_{q_1}[X]$ as well.

Since $(q_1, q_2)=1$ implies that the elements $x \in \Bn_{q_1}[X]/(M_2(X))$ are invariant exactly under $\nu^{q_1q_2}$, using theorem \ref{satz:unteralgebren} we conclude
\begin{equation}\label{eq:iso_spuren}
\Bn_{q_1}[X]/(M_2(X)) \cong \Bn_{q_1q_2} \cong \Bn_{q_1}[X]/(M_1(X)).
 \end{equation}

As $\theta^{(q_1q_2)}=X \mod M_1(X)\Bn_{q_1}[X]$ and $\theta^{(q_2)}=X \mod M_2(X)\Bn_{q_1}[X]$ holds, using isomorphism \eqref{eq:iso_spuren} we see there exists a polynomial $W(X) \in \Bn_{q_1}[X]$ with $\deg2(W(X))<q_2$ such that $W(\theta^{(q_2)})=\theta^{(q_1q_2)}$ holds.\\

Since the isomorphism
\[ 
\alpha_{q_2}: \Bn_{q_2} \rightarrow \Cn_{q_2},\quad \theta^{(q_2)} \mapsto \sum_{i=0}^{q_2-1} a_i \sigma^i(\zeta_\ell^{(q_2)})
 \]
is assumed to be known, we obtain
\[ 
\alpha_{q_1q_2}(\theta^{(q_1q_2)})=\alpha(\theta^{(q_1q_2)})=\alpha(W(\theta^{(q_2)}))=\alpha(W)\left ( \sum_{i=0}^{q_2-1} a_i \sigma^i(\zeta_\ell^{(q_2)}) \right ).
 \]
Since $W(X) \in \Bn_{q_1}[X]$, the coefficients of $W$ depend on $\theta^{(q_1)}$. Using the isomorphism $\alpha_{q_1}$, which is also assumed to be known, we can specify them in terms of $\zeta_\ell$. Finally, this yields $\alpha_{q_1q_2}(\theta^{(q_1q_2)})$ as a function of $\zeta_\ell$.\\

The isomorphism
\[ 
\alpha: \Bn \rightarrow \Cn,\quad \theta \mapsto \sum_{i=0}^{(\ell-1)/2} a_i\zeta_\ell^{c^{2i}}.
 \]
in question is obtained by inductively repeating this procedure.

\subsection{Determination of the trace}
Now we wish to determine the trace $t \mod \ell$ of the Frobenius homomorphism using this isomorphism. Again we set $\theta:=T + (E_P(T))$ in the algebra $\mathbf{B}$, which we use to define $P:=(\theta, \gamma)$, and $\zeta_\ell:=U + (K(U))$ in $\mathbf{C}$. 
Instead of computing $\varphi_p(\theta)$ we now calculate
\begin{equation}\label{eq:spur_best_anfang}
\alpha(\varphi_p(\theta))=\varphi_p(\alpha(\theta))=\varphi_p\left(\sum_{i=1}^{(\ell-1)/2}b_i\zeta_\ell^{c^{2i}}\right)
=\sum_{i=1}^{(\ell-1)/2}b_i^p\zeta_\ell^{pc^{2i}}.
\end{equation}
Obviously, $\zeta_\ell^{pc^{2i}}=\zeta_\ell^k$, where $k \equiv pc^{2i} \mod \ell$ holds. Possibly, the power $\zeta_\ell^k$ with $k<\ell$ still has to be reduced modulo  $K(U)$. Hence, these values can be calculated with negligible cost and the run-time of this step is dominated by the exponentiations of the $b_i$, which lie in $\mathbf{A}$, though. The same applies to the computation of $\alpha(\varphi_p^2(\theta))$. Using the polynomial $G_t$ we can now determine the right hand side of equation \eqref{eq:char_gleichung_x}.\\
In order to determine the left hand side we recall that according to equation \eqref{eq:multi_poly} in addition to the polynomial $G_p$ satisfying $G_p(\theta)=(pP)_x$ there exists another polynomial $H_p$ satisfying $\gamma\cdot H_p(\theta)=(pP)_y$.\\
Using the general formula for adding points on elliptic curves we can now proceed to the determination of the left hand side of equation \eqref{eq:char_gleichung_x}. Assuming that $\phi_p^2(P)\neq \pm pP$ holds (otherwise $t$ can easily be computed) the formula yields
\begin{equation}\label{eq:punkt_addition}
(\phi_p^2(P) + pP)_x=\left(\frac{\varphi_p^2(\gamma)-(pP)_y}{\varphi_p^2(\theta)-(pP)_x} \right)^2 - \varphi_p^2(\theta)-(pP)_x.
\end{equation}
Now
\begin{equation}\label{eq:wurzel_ziehen}
\begin{split}
(\varphi_p^2(\gamma)-(pP)_y)^2&=(\gamma^{p^2}-\gamma\cdot H_p(\theta))^2=\gamma^2(\gamma^{p^2-1}-H_p(\theta))^2\\
&=f(\theta)\left(f(\theta)^{(p^2-1)/2}-H_p(\theta) \right)^2
\end{split}
\end{equation}
holds, so the numerator of the fraction and thus the value $(\phi_p^2(P) + pP)_x$ itself only depend on $\theta$.\\

Having determined the values $\varphi_p^2(\theta), \varphi_p^2(f(\theta))$ using the special structure of $\Cn$ as well as $ G_p(\theta), H_p(\theta)$, we only have to extract one root in \eqref{eq:wurzel_ziehen} to obtain the left hand side of equation \eqref{eq:char_gleichung_x}. However, since this requires non-negligible cost, we slightly modify our approach. Writing
\[ 
A(\theta)=\varphi_p^2(\theta)-G_p(\theta),\ C(\theta)=\varphi_p^2(\theta)+ G_p(\theta)
 \]
and using equations \eqref{eq:punkt_addition} and \eqref{eq:wurzel_ziehen}, we obtain

\begin{align}
& G_t(\varphi_p(\theta))=\frac{\varphi_p^2(f(\theta))-2f(\theta)^{(p^2+1)/2} H_p(\theta) + f(\theta) H_p^2(\theta)}{A^2(\theta)}-C(\theta)\nonumber \\
\Rightarrow &\underbrace{(G_t(\varphi_p(\theta)) + C(\theta))A^2(\theta) - \varphi_p^2(f(\theta)) - f(\theta) H_p^2(\theta)}_{B(t, \theta)} = -2f(\theta)^{(p^2+1)/2} H_p(\theta)\nonumber \\
\Rightarrow & B(t, \theta)^2 = 4f(\theta)^{p^2+1}  H_p^2(\theta)= 4\varphi_p^2(f(\theta))f(\theta)H_p^2(\theta).\label{eq:char_gl_ohne_wurzel}
\end{align}

Obviously, both sides of this equation can be computed without a root extraction and after applying $\alpha$ to \eqref{eq:char_gl_ohne_wurzel} we are left to find the value of $t$ satisfying this equation.\\

Since we only consider $x$-coordinates, we will glean two solutions $\pm t$. In some cases \cite[p. 1251]{Dewaghe} enables us to compute the correct sign of $t$.\\

\section{Run-time}
Our basic operations will be multiplications in $\F_p$. We denote by $\M(n)$ the cost for multiplying two polynomials of degree less than $n$ and by $\Cpoly(n)$ the one for computing $g(h) \mod f$, where $f, g, h \in \F_p[X]$ are of degree less than $n$. Using fast arithmetic one can take (cf. \cite{VzGe, MiMoSc})
\[ 
\M(n)=\tilde{O}(n),\quad \Cpoly(n)=O(n^{(\omega+1)/2}),
 \]
if $n\times n$-matrices over $\F_p$ can be multiplied in $O(n^{\omega})$ operations.

\subsection{Computation of \texorpdfstring{$\beta(\chi_q)$}{beta(chi-q)}}
First, we have to compute the elliptic Gau{\ss} sums. For a character $\chi_q$ of order $q$ we have to determine
\[ 
\tau_e(\chi_q)=\sum_{a=1}^q \rho^a_q \nu^a(\theta^{(q)})=\sum_{a=1}^q \rho_q^a \nu^a \left ( \sum_{j=1}^{(\ell-1)/(2q)} \nu^{jq}(\theta)  \right ).
 \]
We know that the action of $\nu$ on $\theta$ is encoded by the polynomial $G_c(X)$, which can be derived from the division polynomials, where $c$ is a generator of $\F_\ell^*$. First, we wish to compute $\theta^{(q)}$. For this purpose, we need the polynomial $G_{c^q}$ encoding the action of $\nu^q$, determining which requires $O(\log q \Cpoly(\ell))$ operations in $\An$.
Subsequently, we proceed inductively by computing
\[ 
1.\ U:=\theta+G_{c^q}(\theta),\quad 2.\ U:=U+U\circ G_{c^{2q}}=\theta+G_{c^{q}}(\theta)+G_{c^{2q}}(\theta)+G_{c^{3q}}(\theta),\ \ldots.
 \]
The algorithm can be directly adopted from \cite[p.~5]{MiMoSc}. It requires $\log(\frac{\ell-1}{2q})$ steps, each of which has run-time $O(\Cpoly(\ell))$. Afterwards we can determine the value $\tau_e(\chi_q)$ in another $q$ steps of equal run-time. Finally, these calculations require $O(q\Cpoly(\ell))$ operations in $\An$ or 
\begin{equation}
O(q\Cpoly(r\ell))\label{eq:lz1}
\end{equation}
 operations in $\F_p$. The cost for computing the cyclotomic Gau{\ss} sums $\tau(\chi_q)$ is of the same scale since these may be computed similarly.\\
Next, the $q$-th powers of the quantities $\tau_e(\chi_q), \tau(\chi_q)$ have to be computed, which corresponds to $\log q$ multiplications in $\An[\rho_q, \theta^{(q)}]$ and thus 
\begin{equation}
O(\M(rq^2)\log q)\label{eq:lz2}
\end{equation} 
operations. After that, we have to extract one $q$-th root, which according to \cite{DoSc} produces
\begin{equation}
O(\M(q)\M(rq)\log p +q \Cpoly(rq)+\Cpoly(q)\M(rq)\log rq)\label{eq:lz3}
\end{equation}
operations. The value $\beta(\chi_q)$ is now determined. Since $q \in O(p)$, the cost from \eqref{eq:lz2} is dominated by the one from \eqref{eq:lz3}. Since $\beta(\chi_q)$ has to be calculated for all $q \mmid \frac{\ell-1}{2}$, we deduce a total run-time of
\begin{equation}
\sum_{q \mmid \frac{\ell-1}{2}} \eqref{eq:lz1} + \eqref{eq:lz3}
\end{equation}
for this step. Since $\frac{\ell-1}{2}$ has at most $\log \ell$ different prime divisors, we bound the cost by
\begin{align*}
&O(\log \ell (q\Cpoly(r\ell) + \M(q)\M(rq)\log p +q \Cpoly(rq)+\Cpoly(q)\M(rq)\log rq))\\ 
=\ &\tilde{O}(q\Cpoly(r\ell) + rq^2\log p),\label{eq:lz_beta}\numberthis
\end{align*} 
where now $q=\max_i\{q_i \mmid \frac{\ell-1}{2} \}$.
\subsection{Computation of \texorpdfstring{$\alpha$}{alpha}}\label{eq:lz_alpha_2}
Using the approach from section \ref{sec:alpha_best_besser} the cost for determining $\alpha$ may be assessed as follows. First, we determine the isomorphisms $\alpha_q$. Once the value $\beta(\chi_q)$ for $q \mmid \frac{\ell-1}{2}$ is known, one is left to compute the values $\beta(\chi_q^j), j=2, \ldots, q$. This requires $O(q)$ multiplications in $\An[\rho_q, \theta^{(q)}]$ (cf. section \ref{sec:iso_rechnung}). Subsequently, we have to solve a linear system of dimension $q$ over $\An[\rho_q]$, where the matrix to be inverted pertains to a discrete Fourier transform, wherefore this step requires $O(q \log q)$ multiplications in $\An[\rho_q]$. Denoting by $q$ the maximal prime power dividing $n$ we obtain the total cost
\begin{equation}
\tilde{O}(q\M(rq^2)+q\log q \M(rq))=\tilde{O}(rq^3)\label{eq:lz6}.
\end{equation}
For the inductive computation of the isomorphism $\alpha$ from the various $\alpha_q$ the polynomials $M_1, M_2, W$ have to be determined in each step. It should be possible to adopt the algorithm from \cite[pp.~ 4-5]{MiMoSc} to achieve this, which yields the run-time 
\begin{equation*}
O(q_2^{1/2}\M(r\ell)+q_2^{(\omega-1)/2}r\ell),
\end{equation*}
where $\omega$ again denotes an admissible exponent for matrix multiplication, hence $\omega \approx 2.4$ for the asymptotically fastest algorithms \cite{CoWi}. Since $\omega>2$, we bound the run-time by
\begin{equation}
\tilde{O}(q_2^{(\omega-1)/2}r\ell).\label{eq:lz6_2}
\end{equation}
Now we dispose of a polynomial $W=\sum_{k=0}^{q_2-1}w_kX^k \in \Bn_{q_1}[X]$, such that $W(\theta^{(q_2)})=\theta^{(q_1q_2)}$. Since the terms one obtains when computing powers of $\alpha_{q_2}(\theta^{(q_2)})=\sum_{i=0}^{q_2-1} a_i\sigma^i(\zeta_\ell^{(q_2)})$ can in general not directly be written in this form, we use the representation $\alpha_{q_2}(\theta^{(q_2)})=\sum_{i=0}^{\ell-1} a_i\zeta_\ell^i$, which is easily obtained by expanding the $\zeta_\ell^{(q_2)}$, from the start.\\
Since $W$ has degree $q_2-1$, powers have to be computed up to this exponent, which produces $O(q_2\M(r\ell))$ operations. Afterwards, one is left to determine $\sum_{k=0}^{q_2-1} w_k \alpha_{q_2}(\theta^{(q_2)})^k$. The $w_k$ are polynomials in powers of $\alpha_{q_1}(\theta^{(q_1)})$, which may be computed using $O(q_1 \M(r\ell))$ operations by the same reasoning. The final computation of the sum again requires $O(q_2 \M(r\ell))$ operations. Now $\alpha_{q_1q_2}$ is computed. This yields the additional run-time
\begin{equation}
\tilde{O}((q_1+q_2)r\ell),\label{eq:lz7}
\end{equation}
which dominates \eqref{eq:lz6_2}, since $\omega\leq 3$ holds, and is thus the cost for any one of the inductive steps.\\

Obviously, at most $O(\log \ell)$ inductive steps have to be performed until $\alpha$ is determined. Hence, the total run-time is equal to \eqref{eq:lz7}, where $q_1, q_2$ denote the maximal values occurring in the process of the stepwise construction of $ \alpha $, i.~e. the values from the last step when $q_1q_2=\frac{\ell-1}{2}$ holds. In the best case $q_1, q_2$ can be chosen of order $O(\sqrt{\ell})$ resulting in a run-time of $O(r\ell^{3/2})$ operations, whereas the worst-case run-time $O(r\ell^2)$ is attained when the greatest prime power divisor $q \mmid \frac{\ell-1}{2}$ is of order $O(\ell)$.

\subsection{Computation of \texorpdfstring{$t$}{t}}
Due to the special structure of  $\Cn$ the computation of $\varphi_p$ in this algebra can be performed very fast. To be able to efficiently apply $\varphi_p$ to any $b=\sum_{i=0}^{\ell-1} a_i\zeta_\ell^i$ it suffices to precompute the action of $\varphi_p$ on a basis of $\An$. If a power basis $\{1, x, \ldots x^{r-1}\}$ is chosen, it suffices to compute $\varphi_p(x)$ in $\An$, since $\varphi_p(x^k)=(\varphi_p(x))^k$ can be directly derived therefrom. The cost for this precomputation thus amounts to $O(\M(r) \log p)$ operations. Subsequently, the value $\varphi_p(a)=\sum_{j=0}^{r-1} c_j\varphi_p(x)^j$ for $a=\sum_{j=0}^{r-1} c_jx^j \in \An$ can be calculated using $O(r^2)$ multiplications in $\F_p$. To determine $\varphi_p(b)$ this calculation has to performed  $\ell$ times. Hence, the values arising from the action of the Frobenius homomorphism can be computed in run-time $O(r^2\ell)$.\\

Finally, $t$ has to be determined by testing possible values. This requires $O(\ell\M(r\ell))$ operations. The total run-time of this step is thus
\begin{equation}
O(\M(r)\log p + r^2\ell + \ell\M(r\ell))=\tilde{O}(r\log p + r\ell^2).\label{eq:lz8}
\end{equation}

\subsection{Total run-time}
We now combine the run-time estimates from the preceding sections, which yields
\begin{align*}
&\tilde{O}(q\Cpoly(r\ell) + rq^2\log p + rq^3 + (q_1+q_2)r\ell + r\log p + r\ell^2)\\
=\ &\tilde{O}(rq^2\log p + q\Cpoly(r\ell) +rq^3 + r\ell^2)\numberthis \label{eq:lz_gesamt_2}.
\end{align*}

We compare this to the run-time for the Elkies case, which amounts to $\tilde{O}(\ell \log p)$ operations, whereas the corresponding computations in Schoof's original algorithm required $\tilde{O}(\ell^2 \log p)$ operations. Since we expect the largest prime $\ell$ to be considered to be of size $O(\log p)$, the original algorithm allows to handle Atkin primes $\ell \in \tilde{O}(\sqrt{\log p})$. We now consider the conditions that arise if we require that our run-time not exceed $\tilde{O}(\log^2 p)$, i.~e. the cost for Elkies primes  $\ell \in O(\log p)$.\\

This yields the following restrictions:
\begin{align*}
 rq^2 \in\ &\tilde{O}(\log p),\\
 q\Cpoly(r\ell) \in\ &\tilde{O}(\log^2 p),\quad rq^3 \in \tilde{O}(\log^2 p),\quad  r\ell^2 \in \tilde{O}(\log^2 p).
\end{align*}

The first condition implies the middle one in the second line. Assuming $r$ to be small, one obtains
\[ 
q  \in \tilde{O}(\sqrt{\log p}),\quad q\Cpoly(\ell) \in \tilde{O}(\log^2 p),
\]
so provided $q$ is small enough one might use Atkin primes satisfying $\Cpoly(\ell)\in \tilde{O}((\log p)^{1.5})$, i.~e., $\ell\in \tilde{O}((\log p)^{0.89})$ when using the best known $\omega \approx 2.4$. However, if $r$ is small, one might as well use the existing generic method for Atkin primes (cf. \cite{Schoof, Mueller}), which is efficient in this case.\\
Abandoning the assumption that $r$ is small and thus assuming only $r \in O(\ell)$, we observe that on the one hand $q$ has to take smaller values and on the other hand $O(q\Cpoly(r\ell))$ becomes the dominating term in the second line. In case $q$ is small, the bound for $\ell$ amounts to
\[\ell \in \tilde{O}((\log p)^{\frac{2}{\omega+1}}), \]
which leads to $\ell \in \tilde{O}\left((\log p)^{0.59}\right)$ instead of $\ell \in \tilde{O}\left((\log p)^{0.5}\right)$ from Schoof's algorithm. We conclude the method might allow to use slightly larger values for $\ell$ if the corresponding $r$ and $q$ are not too large. For the time being, this gain remains theoretical, though, since the method has not been implemented for large-scale computations.\\

\subsection*{Acknowledgements}
I would like to thank Jean-Pierre Flori for a helpful remark concerning an earlier version of this paper.

\footnotesize{
\bibliography{lit}}

\end{document}